\documentclass{article}
\usepackage{amsthm}
\usepackage{amsmath}
\usepackage{amssymb}

\newtheorem{theorem}{Theorem}[section]

\theoremstyle{definition}

\newtheorem{example}[theorem]{Example}

\theoremstyle{remark}
\newtheorem{remark}[theorem]{Remark}

\numberwithin{equation}{section}

\newtheorem{definitons}{Definitions}
\def\O{{\mathcal O}}
\def\Ac{{\mathcal A}}
\def\Bc{{\mathcal B}}
\def\Mc{{\mathcal M}}

\def\Oc{{\mathcal O}}

\def\P{{\mathbb P}}

\def\Hom{{\mathcal H om}}
\def\banica{{B\u anic\u a}}

\def\qed{\hfill$\Box$\vskip10pt}

\begin{document}

\title{A Class of Locally Complete Intersection Multiple Structures on Smooth
Algebraic Varieties as Support}

\author{Nicolae Manolache}
\maketitle
\noindent\hfil {Math. Sub. Class.: Primary 14M05, 13H10}

\date{}

\begin{abstract}
The systematic study of multiple structures, made necessary (and possible) by the notion of {\bf scheme}
introduced by Grothendieck, began with the papers 
\cite{Fo}, \cite{Fe} and was continuated in \cite{BF1} and \cite{M1}, \cite{M2}, \cite{M3}
(the last three use also ideeas from \cite{PS}). In \cite{BF1}, \cite{BF2}, along the 
classification, up to multiplicity $4$,
of multiple locally complete intersection (lci for short) structures 
on a smooth curve embedded in a smooth threefold, general classes of multiple 
structures are introduced, the so-called ``primitive'' and ``quasiprimitive'' 
structures. The primitive ones are characterized by the fact that, locally, they 
are defined by equations of the type $x^n=0$, $y=0, z=0, \ldots ,u=0$; 
the quasiprimitive ones are those which are generically primitive. 
The general study of these structures was continued by several authors, from 
which we mention  \cite{Bo}, \cite{Dr1}, \cite{Dr2}, \cite{Dr3}, \cite{Dr4}.
 
In this paper we give the construction of a class of multiple locally complete 
intersection (lci for short) structures on a smooth algebraic variety as 
support. This class contains the lci structures defined locally by 
the ``the next`` monomial equations, namely those which are defined locally 
by  equations of the form $x^n=0$, $y^2=0, z=0,\ldots ,u=0$.
\end{abstract}

\maketitle

\section{Preliminaries}

Let $X$ be a smooth connected algebraic variety over an algebraically closed 
field $k$.
A (locally) Cohen-Macaulay scheme $Y$ is called a \emph{multiple structure on $X$} 
if the subjacent 
reduced scheme $Y_{red}$ is $X$. In this case all the local rings of $Y$ have the 
same multiplicity
(cf. \cite{M2}), which is called \emph{the multiplicity of $Y$}. Let $Y$ be 
embedded in a smooth variety $P$. 
Let $I$ be the (sheaf) ideal of $X$ in $P$ and $J$ be the ideal of $Y$ in $P$. 
Let $m$ be the positive integer such that $I^m \not\subset J$,
$I^{m+1} \subset J$. To  $Y$ one associates canonically three filtrations. They are:\\
1. Let $I^{(\ell)}$ be the ideal obtained throwing away the embedded components 
of $I^\ell+J$ and let 
$Z_\ell$ be the corresponding scheme. This  gives the 
\emph{\banica -Forster filtration} (cf. \cite{BF2}):
\begin{displaymath}
\begin{array}{ccccccccccccc}
 \O_Y=I^{(0)} & \supset & I= I^{(1)} & \supset &  I^{(2)} & \supset &  \ldots &  
\supset &   I^{(m)} & 
\supset & I^{(m+1)}=0\medskip \\
&  & X=Z_1 & \subset & Z_2 & \subset &  \ldots &  \subset & Z_m & \subset & 
Z_{m+1} =  Y
\end{array}
\end{displaymath}
$Z_\ell$ are not, in general, Cohen-Macaulay. But this is true if ${\rm dim}(X)=1$.
The graded associated object $\Bc (Y)=\bigoplus_{\ell =0}^m I^{(\ell)} / I^{(\ell +1)} $ 
is naturally a 
graded $\O_X$-algebra. If the schemes $Z_\ell$ are Cohen-Macaulay, the graded components 
of $\Bc (Y)$ are locally 
free sheaves on $X$.\\
2. Let $X_\ell $ be defined by $I_\ell =J:I^{m+1-\ell}$.  When $X_\ell$ are Cohen-Macaulay, 
the quotients 
$I_\ell /I_{\ell +1}$ are 
locally free sheaves on $X$. Again, if ${\rm dim}(X)=1$, 
$X_\ell $ are Cohen-Macaulay.
This is also true if $Y$ is lci (i.e. locally complete intersection) of multiplicity
at most $6$ (cf. \cite{M2}).
In general this is not always the case. This filtration was considered in \cite{M1}.\\
3. Let $Y_\ell$ be the scheme given by $J_\ell=J:I_{m+1-\ell}=J:(J:I^\ell)$. When $X_\ell$ 
is Cohen-Macaulay, $Y_\ell$
has the same property. The graded object 
$\Ac (Y)=\bigoplus_{\ell =0}^m J_{\ell} / J_{\ell +1}$ is a graded 
$\O _X$-algebra and $\Mc (Y) =\bigoplus_{\ell =0}^m I_\ell / I_{\ell+1}$ is a 
graded $\Ac(Y)$-module.
This filtration was considered in \cite{M2}.\\
The system of the graded components ($\Ac_0(Y),\ldots  \Ac_m(Y);\Mc_0(Y),\ldots  \Mc_m(Y))$ 
is called 
\emph{the type of $Y$}. $Y$ is called \emph{of free type} when all the graded pieces are 
locally free.
As already remarked, in dimension $1$, or if $Y$ is lci of multiplicity up to $6$, 
this is  the case.

The previous filtrations have the following properties:\\
1) In general the above filtrations are different.
Take for instance $X=Spec(k)$, $Y=\hbox{Spec}(k[x,y]/(x^3,xy,y^4))$, $P=\hbox{Spec}(k[x,y])$\\
2) $Z_\ell \subset Y_\ell \subset X_\ell$ \\
2') there are canonical morphisms: $\Bc(Y)\to \Ac(Y)\to \Mc (Y)$\\
3) The multiplications
\begin{displaymath}
\begin{array}{lll}
\Ac_{\ell_1}(Y)\otimes \Ac_{\ell_2}(Y) & \to & \Ac_{\ell_1+\ell_2}(Y)\medskip\\
\Ac_{\ell_1}(Y)\otimes \Mc_{\ell_2}(Y) & \to & \Mc_{\ell_1+\ell_2}(Y)
\end{array}
\end{displaymath}
are never the zero maps for $\ell_1, \ell_2 \ge 0$, $\ell_1+\ell_2 \le m$ (cf. \cite{M2}).\\
4) From the definitions, one has the exact sequences:
\begin{displaymath}
\begin{array}{c}
0\to \Mc_\ell(Y)\to \O_{X_{\ell+1}} \to \O_{X_\ell} \to 0\medskip\\
0\to \Ac_\ell(Y)\to \O_{Y_{\ell+1}} \to \O_{Y_\ell} \to 0
\end{array}
\end{displaymath}
5) If $Y$ is Gorenstein of free type, then $X_\ell $ and $Y_{m+1-\ell}$ are locally 
algebraically linked
(cf. \cite{M1}). In particular one has the exact sequences:
\begin{displaymath}
\begin{array}{c}
0\to \omega_{X_{m+1-\ell}}\otimes \omega_Y^{-1}\to \O_Y \to \O_{Y_\ell} \to 0\medskip\\
0\to \omega_{Y_{m+1-\ell}}\otimes \omega_Y^{-1}\to \O_Y \to \O_{X_\ell} \to 0
\end{array}
\end{displaymath}
6) Let $Y$ be a free type Cohen-Macaulay multiple structure on a smooth support $X$.\\
Then (cf \cite{M3}) $Y$ is Gorenstein if and only if the following conditions are fulfilled :\\
(a) $\Ac_m(Y)$ and $\Mc_m(Y)$ are line bundles\\
(b) $\Ac_m(Y)=\Mc_m(Y)$  \\
(c) The canonical maps:
\begin{displaymath}
\Ac_\ell (Y)\to \Hom _{\O_X}(\Mc_{m-\ell}(Y),\Mc_m(Y))\cong \Mc_{m-\ell}(Y)^\vee\otimes \Mc_m(Y)
\end{displaymath}
are isomorphisms. \\
6') In particular: if $Y$ is Gorenstein of free type, then (cf also \cite{M2}):
\begin{displaymath}
\begin{array}{l}
(a)\  \rm{rank}\ \Ac_\ell(Y)=\rm{rank}\ \Mc_{m-\ell}(Y)\medskip\\
(b)\  \Ac_\ell(Y)=\Mc _\ell (Y) \rm{\ iff\ } \rm{rank}\ \Ac_\ell(Y)=\rm{rank}\ \Ac_{m-\ell}(Y)
\end{array}
\end{displaymath}
In this paper all the schemes are algebraic schemes over a fixed algebraically closed field $k$, 
of characteristic $0$.

\section{Some structures of multiplicity 2n}

\begin{definitons}
1. Let $E$ be a rank $2$ vector bundle on $X$ and $L$ be a line bundle on $X$. If $s$ is 
a section of 
$L\otimes S^2E$, locally of the form $\ell\otimes (ae_1^2+2be_1e_2+ce_2^2)$,
($\ell $ a local generator of $L$ and $e_1$, $e_2$ local generators for $E$), then 
$\delta (s) \in \Gamma
(X,L^2\otimes(det E)^2)$ is the section defined locally by 
$(b^2-ac)\ell ^2\otimes (e_1\wedge e_2)^2$.\\
When $s\in \Hom (L, S^2 E)$ define 
$\delta (s)=\delta (s\otimes 1_{L^{-1}})\in \Gamma (X, L^{-2}\otimes(det E)^2)$.\\
When $s \in \Hom (S^2 E, L)$, define 
$\delta (s) = \delta (s ^\vee )\in \Gamma (X, L^2\otimes (detE)^{-2})$.\\
Extend this definition also to a surjection $\varphi : S^2 E \to F$, where $F$ 
is a vector bundle of rank $2$.
Namely, if $s_\varphi : L\to S^2E$ is the kernel of $\varphi$, then put 
$\delta(\varphi)=\delta (s_\varphi )$.

2. With $L$ and $E$ as above, for $s\in \Gamma (X,L\otimes S^nE)$,
define the Hessian $h(s)\in \Gamma (X, L^2\otimes (det E)^2\otimes S^{2n-4}E)$ in 
the obvious way: if $s$ is given 
locally by $\ell \otimes \sum {n \choose i} a_{ij}e_1^ie_2^j$, then $h(s)$ is locally 
the Hessian of this symmetric form.\\
When $s \in Hom(L, S^n E)$ define 
$h(s)=h(s\otimes 1_{L^{-1}}) \in \Gamma (X, L^{-2}\otimes (det E)^2\otimes S^{2n-4}E)$\\
When $s \in Hom(S^n E, L)$ define 
$h(s)=h(s ^\vee) \in \Gamma (X, L^2\otimes (det E)^{-2}\otimes S^{2n-4}E^\vee)$
\end{definitons}
\begin{theorem}
Let $X \subset P$ be smooth varieties and let $I$ be the ideal sheaf of $X$ in $P$. 
The following construction produces a 
lci multiple
structure $Y$ with support $X$, of multiplicity $2n$.\\
Step 1. Take  a rank $2$ vector bundle $E_1$ on $X$, a surjection $p_1:I/I^2 \to  E_1$ 
and define
$I_2:=ker(I\to I/I^2 \to  E_1)$. Then it follows:
\begin{displaymath}
S^2E_1 \cong I^2/II_2,\ \  S^3E_1 \cong I^3/I^2I_2 \ \ \hbox{etc.}\\
\end{displaymath}
Step 2. Take $E_2$ a rank $2$ vector bundle on $X$ and $p_2: I_2/II_2 \to E_2$, 
such that the map
$\mu _2: S^2 E_1 =I^2/II_2 \hookrightarrow I_2/II_2 \to E_2$ is surjective. 
Take $I_3=ker (I_2 \to 
I_2/II_2 \to E_2)$.Then $E_3:=II_2/II_3$ is a rank $2$ vector bundle on $X$.\\
Step 3. Take  a surjection $p_3:I_3/II_3 \to E_3$ such that the composition of $p_3$ 
with the natural 
inclusion $i_3:II_2/II_3 \to I_3/II_3$ 
is an isomorphism of vector bundles. Define $I_4=ker (I_3 \to I_3/II_3 \to E_3)$. 
Then $E_4:=II_3/II_4$ is a rank $2$ vector 
bundle on $X$.\\ 
\smallskip \\
\vdots \medskip \\ 
Step k. ($4\le k\le n-2$) Take  a surjection $p_k:I_k/II_k \to E_k$ such that the composition of $p_k$ with the natural 
inclusion $i_k:E_k = II_{k-1}/II_k \to I_k/II_k$ 
is an isomorphism of vector bundles. Define $I_{k+1}=ker (I_k \to I_k/II_k \to E_k)$. 
Then $E_{k+1}:=II_k/II_{k+1}$ is a rank $2$ vector 
bundle on $X$.\medskip \\ 
\vdots \medskip \\ 
Step n-1. Take  a surjection $p_{n-1}:I_{n-1}/II_{n-1} \to E_{n-1}$ such that the composition of $p_{n-1}$
 with the natural inclusion $i_{n-1}:E_{n-1} = II_{n-2}/II_{n-1} \to I_{n-1}/II_{n-1}$ 
is an isomorphism of vector bundles. Define $I_n=ker (I_{n-1} \to I_{n-1}/II_{n-1} \to E_{n-1})$. 
\medskip \\ 
Step n. Take  a line bundle $L$ on $X$ such that 
$c_{2n-3}(L^2\otimes (detE_1)^{-2}\otimes S^{2n-4}E^\vee _1)=0 $.
Take a surjection $p_n: I_n/II_n\to L$
such that $\mu _n :=p_nj_n$ is surjective, where $j_n$ is the canonical map 
$j_n:S^nE_1 \to I_n/II_n$ and such that the 
Hessian of $\mu _n$ vanishes nowhere. Define $I_{n+1}:=ker(I_n \to I/II_n \to L)$. 
Then $I_{n+1}$ defines a subscheme 
$Y$ of $P$ which is lci of multiplicity $2n$.\\

One has: $\omega _Y \cong \omega _X \otimes L^{-1}$ and $\Bc (Y)=\Ac (Y) =\Mc (Y) = 
\Oc _X \oplus E_1 \oplus \ldots \oplus E_{n-1}\oplus L$. Moreover, by \cite{M3}, the 
canonical maps
$E_\ell \to \Hom (E_{n-\ell}, L)$ are isomorphisms.
\end{theorem}
\emph{Proof.}
We will determine at each step the local equations of the scheme defined by 
$I_\ell$, i.e. 
the equations in the {\em completions} of the local rings with respect to their 
maximal ideals.

\emph{Step 1.} It is clear that in each point of $X$ there are local parameters 
$x,y,z,\ldots u$ 
which define the scheme $X$ in that point such that the ideal $I_2$ is of the 
form $(x^2,xy,y^2, z,\ldots u)$.
As $E_1=I/I_2$, the multiplication map 
\begin{displaymath}
 E_1\otimes E_1 \cong {I \over I_2} \otimes {I \over I_2} \to {I^2 \over II_2}
\end{displaymath}
factors through $E_1\otimes E_1 \to S^2 E_1$ and thus the map $S^2 E_1 \to I^2/II_2$, 
as a surjection between
rank $3$ vector bundles is an isomorphism. In the same way one shows all the 
other isomorphisms.

\emph{Step 2.} A local computation shows that the following three situations can occur:\\
a) $\delta (\mu _2)=0$. Then, in convenient local coordinates, 
$I_3=(x^3, x^2y,y^2,z,\ldots ,u)$.\\
a') $\delta (\mu _2)$ defines an effective nonzero divisor $D$. Then in the points of this 
divisor, in convenient coordinates,
$I_3=(x^3, x^2y,y^2+\alpha _2 x^2,z,\ldots ,u)$, where $\alpha _2=0$ 
is the local equation of $D$.\\
b) $\delta (\mu _2)$ vanishes nowhere. After a change of variables, one has \\
$I_3=(x^3,xy,y^3,z,\ldots u)$.\\
In all these situations, by local computations, one shows that $E_3:=II_2/II_3$ is a 
rank $2$ vector bundle.

\emph{Step 3.} We treat in the same time the situations a) and a') from step 2.
One shows easily that, in new local coordinates (i.e. new generators of the maximal
ideal):
$I_4=(x^4, x^3y,y^2+\alpha _2 x^2+\alpha _3 x^3,z,\ldots ,u)$
In the case b), one obtains $I_4=(x^4,y^4, xy-\lambda x^3-\mu y^3, z,\ldots ,u)$. 
Via the change of coordinates
$x=X+\mu Y^2$, $y=Y+\lambda X^2$, one gets 
$I_4=(X^4,Y^4, XY, z, \ldots ,u)$.
In both cases, by local computations, one shows that $E_4:=II_3/II_4$ is a 
rank $2$ vector bundle.\\
\vdots

\emph{Step $k$. $(4\le k \le n-2)$} 
By induction, at steps $k-1$ and $k-2$ we get two possibilities for the local shape of
our ideals:

a) $I_k=(x^k,x^{k-1}y, y^2+\alpha _2x^2+\ldots +\alpha _{k-1}x^{k-1},z, \ldots ,u)$,

\ \ \ $I_{k-1}=(x^{k-1},x^{k-2}y, y^2+\alpha _2x^2+\ldots +\alpha _{k-2}x^{k-2},z, \ldots ,u)$

or

b) $I_k=(x^k,xy,y^k,z,\ldots ,u)$,

\ \ \ $I_{k-1}=(x^{k-1},xy,y^{k-1},z,\ldots ,u)$,

In case a), as the local generators of $E_k=II_{k-1}/II_k$ are $x^k$, $x^{k-1}y$, one obtains,
in convenient coordinates,  $I_{k+1}=(x^{k+1},x^ky, y^2+\alpha _2x^2+\ldots +\alpha _kx^k,z, \ldots ,u)$

In case b), as the local generators of $E_k=II_{k-1}/II_k$ are $x^k$, $y^k$ one obtains
$I_{k+1}=(x^{k+1},xy-\lambda x^k-\mu y^k ,y^{k+1},x^2y,xy^2,z,\ldots ,u)$. 
With the change $y=Y+\lambda x^{k-1}$ one gets $I_{k+1}=(x^{k+1}, xY-\mu (Y+\lambda x^{k-1})^k ,
(Y+\lambda x^{k-1})^{k+1}, x^2Y, xY^2, z,\ldots ,u)$, where $z$, \ldots $u$ are changed tacitly.
Now, taking $x=X+\mu Y^{k-1}$, one gets $I_{k+1}=(X^{k+1},XY,Y^{k+1},z,\ldots ,u)$

In both cases, by local computations, one shows that $E_{k+1}:=II_{k-1}/II_{k}$ is a 
rank $2$ vector bundle on $X$.\\

\emph{Step $n-1$} is similar with the previous one, but we are no more interested in 
the vector bundle $E_n:=II_{n-1}/II_n$.

\emph{Step $n$}. We have two cases:

a) $I_n=(x^n, x^{n-1}y, y^2+\alpha _2 x^2 +\ldots +\alpha _{n-1}x^{n-1}, z, \ldots ,u)$,

or 

b) $I_{n}=(x^n,xy,y^n,z,\ldots ,u)$.

Take $\ell $ to be the local generator of the line bundle $L$.

In case a), $\mu _n$ will have the shape:

\begin{displaymath}
 \mu_n(x^{n-2p}y^{2p})=(-\alpha _2)^pr\ell, \ \ 
\mu_n(x^{n-2p-1}y^{2p+1})=(-\alpha _2)^ps\ell
\end{displaymath}
A direct computation of the Hessian of $\mu_ n$ gives:
\begin{displaymath}
h(\mu_n)\equiv -(n-1)^2s^2\ell ^2x^{2n-4}\ \hbox{mod } \alpha _2
\end{displaymath}
From here it follows that the condition upon Hessian implies that $s$ is invertible 
and then, after some change of coordinates:
\begin{displaymath}
I_{n+1}=(x^n, y^2+\alpha _2x^2+\ldots +\alpha _{n-1}x^{n-1},z,\ldots ,u)
\end{displaymath}

b) In this case the application $\mu _n$ has the form:
\begin{displaymath}
\mu _n(x^n)=r\ell, \ \ \mu _n(y^n)=s\ell , \ \ \mu _n(x^{n-p}y^p)=0 \hbox{ for } 
p\neq 0 \hbox { or } p\neq n
\end{displaymath}
Then the Hessian will be:
\begin{displaymath}
h(\mu_n)=n^2(n-1)^2rs\ell ^2x^{n-2}y^{n-2}\ ,
\end{displaymath}
which shows that both $r$ and $s$ should be  invertible.

Then, after a change of coordinates:
\begin{displaymath}
I_{n+1}=(x^n+y^n,xy,z,\ldots ,u)
\end{displaymath}
\qed

\begin{remark}
Conversely, all  multiple structures $Y$ on a smooth variety $X$ embedded in a smooth variety $P$, 
given locally by
ideals of the form $J=(x^n, y^2+\alpha _2x^2+\ldots +\alpha _{n-1}x^{n-1},z,\ldots ,u)$,
or $J=(x^n+y^n,xy,z,\ldots ,u)$ are obtainable by the procedure described in the above theorem.
\end{remark}

\begin{remark}
 The case $n=2$ for curves in threefolds is treated  in \cite{BF1}, \cite{BF2}, 
the case  $n=3$ is covered in \cite{M2}.
\end{remark}

\begin{example}
The above construction has almost no degree of liberty in codimension $2$.
We shall exemplify this for the case of a line embedded in $\P^3 $.
Take $X$ to be a line defined by the homogeneous ideal $I=(x,y)$ and consider $u,v$ 
to be the homogeneous 
coordinates on $X$. 
We shall describe the above construction for $n=2$ (whence multiplicity $4$) and $n=3$ 
(which gives a multiplicity $6$ structure).

\emph{Case $n=2$}.
$E_1$ has to be $I/I^2=2\O_X (-1)$ and then $I_2=(x^2,xy,y^2)$. The condition on $L$ is 
fullfilled only 
by $L=\O_X(-1)$. Then $I_3=(x^2,y^2)$, in  new coordinates.

\emph{Case $n=3$}. $E_1$ and $I_2$ are as above. The condition on the Hessian of $\mu _3$ 
determines the line bundle $L$:
$L=\O _X(-3)$. Then $E_2\cong E_1^\vee \otimes L =2\O_X(-2)$.
The ideal $I_3$ will have the shape $I_3=(c_1x^2+c_2xy+c_3y^2, x^3,x^2u,xy^2,y^3)$, 
where $c_1,c_2,c_3$ are constants.
After a change of the homogeneous coordinates in $\P^3$ this ideal takes one of 
the shapes:
$I_3=(x^2,xy^2,y^3)$ or $I_3=(x^3,xy,y^3)$.
Accordingly, for $I_4$ there are two possibilities:
$I_4=(x^2,y^3)$ or $I_4=(xy, x^3+y^3)$.
\end{example}

In the case of codimension greater than $2$ one can give examples which are not globally 
complete intersections. The simplest are with $n=2$.
\begin{example}
Consider a line $X$ embedded in $\P^4$. Then:
\begin{displaymath}
 {I \over I^2}={(x,y,z) \over (x,y,z)^2}=3\O_X(-1)
\end{displaymath}
We take $E_1=\O_X(-1)\oplus \O_X(r)$, with $r\ge 0$. Consider $p_1: I/I^2 \to E_1$ given 
by a matrix of the shape
\begin{displaymath}
\begin{pmatrix} 
 \ 0 & 0 & 1 \\ -b & a & 0
\end{pmatrix} \ ,
\end{displaymath}
where $a$,$b$ are forms of degree $r+1$ in the homogeneous coordinates $u$,$v$ on $X$,
without common zeroes on $X$.
Then $I_2=(ax+by,x^2,xy,xz,yz,y^2,z^2)$. The condition on the Hessian of $\mu _2$ determines
$L=\O_X(r-1)$. We take $p_2$ such that the map 
$\mu _2 :S^2E_1=\O_X(-2)\oplus \O_X(r-1)\oplus \O _X(2r)\to L$
to be the projection on the second factor.
Then $I_3=(ax+by,x^2,xy,y^2,z^2)$
\end{example}
A similar example can be given for the case $dimX=2$:
\begin{example}
 In $\P ^6$ let $X$ be the plane given by the eqations $x=y=z=t=0$.
If $u,v,w$ are homogeneous coordinates on $X$ and $a,b,c$ are forms
in $u,v,w$ of the same degree $r+1$, then $Y$ defined by the ideal
$I_3=(bx-ay, cy-bz, az-cx)+(x,y,z)^2+(t^2)$
is a multiple structure with $n=2$ on $X$.
\end{example}

{\bf Acknowledgment} The author was partially supported by the CEEX Programme of the 
Romanian Ministry of Education and
Research, contract 2-CEx06-11-20/2006 and the Mittag-Leffler Institute. The author 
wants to express his warmest thanks
for the hospitality and the wonderfull working atmosphere and conditions offered by 
the Mittag-Leffler Institute.

\bibliographystyle{amsplain}

\noindent Nicolae Manolache\\
Institute of Mathematics "Simion Stoilow"\\
of the Romanian Academy \\
P.O.Box 1-764
Bucharest, RO-014700

\noindent e-mail: nicolae.manolache@imar.ro
\end{document}